\documentclass[11pt]{article}

\usepackage{enumerate}
\usepackage{amsmath,amsfonts,amssymb,amsthm,graphics,amscd,comment}

\newtheorem{theorem}{Theorem}[section]
\newtheorem{lemma}{Lemma}[section]
\newtheorem{corollary}{Corollary}[section]

\newtheorem{remark}{Remark}[section]

\theoremstyle{definition}
\newtheorem{definition}{Definition}[section]
\theoremstyle{example}
\newtheorem{example}{Example}[section]

\def\summ{\sum\limits}

\begin{document}

\title{Submajorization on $\ell^p(I)^+$ determined by increasable doubly substochastic operators and its linear preservers}

\author{Martin Z. Ljubenovi\' c, Dragan S. Raki\' c\thanks{This research was financially supported by the Ministry of Education, Science
and Technological Development of the Republic of Serbia (Contract No.451-03-9/2021-14/200109) and by the bilateral
project between Serbia and Slovenia (Generalized inverses, operator
equations and applications, Grant No. 337-00-21/2020-09/32)}}

\date{}

\maketitle

\begin{abstract}

We note that the well-known result of Von Neumann \cite{von} is not valid for all doubly substochastic operators on discrete Lebesgue spaces $\ell^p(I)$, $p\in[1,\infty)$. This fact lead us to distinguish two classes of these operators. Precisely, the class of increasable doubly substochastic operators on  $\ell^p(I)$  is isolated with the property that an analogue  of the Von Neumann result on  operators in this class is true. The submajorization relation $\prec_s$ on the positive cone   $\ell^p(I)^+$, when $p\in[1,\infty)$, is introduced by increasable substochastic operator and it is provided that submajorization may be considered as a partial order.  Two different shapes of  linear preservers of submajorization  $\prec_s$ on $\ell^1(I)^+$ and on $\ell^p(I)^+$, when $I$ is an infinite set, are presented. 

\medskip


\end{abstract}


\section{Introduction}
A square  $n\times n$ non-negative   matrix $A=(a_{ij})$ is called a doubly stochastic if all of its row sums and all of its column sums are equal $1$. A square  $n\times n$ matrix $D=(d_{ij})$ with non-negative entries is called doubly substochastic if there is a doubly stochastic matrix $A=(a_{ij})$ such that  
\begin{eqnarray}\label{doublysubstochasticmatrix}
d_{ij}\leq a_{ij},\;\;\; 1\leq i\leq n, \;\;\;  1\leq j\leq n.
\end{eqnarray}
As a direct consequence, it follows that all  row sums and column sums of matrix $D$ are less than or equal to $1$.
 Von Neumann \cite{von} provides the converse in the next theorem.
 
\begin{theorem}\label{vonnn}\cite[Theorem I.2.C.1]{marshall}
For every $n\times n$ non-negative  matrix $D=(d_{ij})$ with 
\begin{eqnarray}\label{doublysubstochasticmatrix2}
\sum_{i=1}^nd_{ij}\leq 1, \;\;\; \forall j\in\{1,2,\ldots,n\}  \;\;\;\text{and}\;\;\sum_{j=1}^nd_{ij}\leq 1,\;\;\; \forall i\in\{1,2,\ldots,n\},  \;\;\;
\end{eqnarray}
 there exists a doubly stochastic matrix $A=(a_{ij})$ such that $d_{ij}\leq a_{ij}$, for all $i,j\in\{1,2,\ldots,n\}$.
\end{theorem}
Thus, the statement \eqref{doublysubstochasticmatrix2}  may be considered as an alternative definition of a doubly substochastic matrix.
However, the last result of von Neumann may not be generalized to  infinite matrices considered as bounded linear operators on discrete Lebesgue spaces $\ell^p(I)$, whenever $p\in[1,\infty)$.
\begin{definition}\cite[Definition 3.1]{martin2}\label{defsub}
 Let $p\in [1,\infty)$ and let $A:\ell^p(J)\longrightarrow
\ell^p(I)$ be a bounded linear operator, where $I$ and $J$ are two
non-empty sets. The operator $A$ is called
\begin{itemize}
    \item[i)] {\em row substochastic}, if $A$ is positive and $\forall i\in I$,  $\sum_{j\in J}\langle
    Ae_j,e_i\rangle\leq 1$.
    \item[ii)] {\em column substochastic}, if $A$ is positive and $\forall j\in J$,  $\sum_{i\in I}\langle
    Ae_j,e_i\rangle\leq 1$.
    \item[iii)] {\em doubly substochastic}, if $A$ is both row and column
    substochastic.
\end{itemize}
\end{definition}

\begin{definition}\cite[\text{Definition 2.1}]{bahrami}\label{DsS}
 Let $p\in [1,\infty)$ and let $A:\ell^p(J)\longrightarrow
\ell^p(I)$ be a positive bounded linear operator, where $I$ and $J$ are
two non-empty sets. The operator $A$ is called
{\em doubly stochastic}, if $$\forall i\in I,\;\;\sum_{j\in J}\langle
    Ae_j,e_i\rangle= 1,\;\;\text{and}\;\;\forall j\in J,\;\;\sum_{i\in I}\langle
    Ae_j,e_i\rangle= 1.$$
\end{definition}

It is easy to see that the left and the right shift operators $L,R:\ell^p(\mathbb{N})\rightarrow\ell^p(\mathbb{N})$ defined by
\begin{eqnarray*}
\langle Le_j,e_i\rangle:= \left\{\begin{matrix}
1,\hspace{0.2cm} j-i=1, \\
 0, \hspace{0.2cm} \text{otherwise}
\end{matrix}\right.\;\;\;\;\;
\langle Re_j,e_i\rangle:= \left\{\begin{matrix}
1,\hspace{0.2cm} i-j=1, \\
 0, \hspace{0.2cm} \text{otherwise}
\end{matrix}\right.
\end{eqnarray*}

 with matrix forms
\begin{eqnarray}\label{LR}
 L=\left[\begin{matrix}
0  &  1 & 0&0&\ldots  \\
0  &  0 & 1&0&\ldots  \\
0  &   0 &  0 &1   &\ldots \\
0  &   0 &  0 &0   &\ldots \\
0  &   0 &  0 &0   &\ldots \\
\vdots & \vdots& \vdots&\vdots  & \ddots  \\
\end{matrix}\right]
\;\;\;\;\;\; 
R=\left[\begin{matrix}
0  &  0 & 0&0&\ldots  \\
1  &  0 & 0&0&\ldots  \\
0  &   1 &  0 &0   &\ldots \\
0  &   0 &  1 &0   &\ldots \\
0  &   0 &  0 &1   &\ldots \\
\vdots & \vdots& \vdots&\vdots  & \ddots  \\
\end{matrix}\right]
\end{eqnarray}
are doubly substochastic, but there is no
 doubly stochastic operator $A$ such that $$\langle Le_j,e_i\rangle\leq \langle Ae_j,e_i\rangle,\;\;\forall i,j\in \mathbb{N}$$ or $$\langle Re_j,e_i\rangle\leq \langle Ae_j,e_i\rangle,\;\;\forall i,j\in \mathbb{N}$$ holds (the definition \eqref{LR} of the matrix operators $L$ and $R$ is supported by Theorem \ref{forallp}).
We conclude that  shift operators $L$ and $R$ are not doubly substochastic in the sense of definition \eqref{doublysubstochasticmatrix}. This fact lead us to consider a subclass of doubly substochastic operators which is called increasable doubly substochastic operators (introduced in Definition \ref{defsubbystochastic}) based on the definition  \eqref{doublysubstochasticmatrix} to get that each operator from this class satisfies conditions \eqref{doublysubstochasticmatrix} and \eqref{doublysubstochasticmatrix2} (reformulated for operators on $\ell^p(I)$).

In recent years, there is a big progress towards developing various extensions of the most important majorization relations  on sequence spaces \cite{arveson,kaftal} and on descrete  Lebesgue spaces \cite{bahrami,bahrami4,bahrami6,bahrami7,bahrami10,martin,martin2,
martin3,martin4,martin5,martin7} with apropriate generalizations of some famous theorems in linear algebra \cite{antezana,antezana2,argerami,kennedy,loreaux,massey,neumann,per}. 
There are a lot of applications of the majorization theory in various branches of mathematics  and  there exist significant conections  with the other  science like physics, quantum mechanics and quantum information theory \cite{nature,li,manjegani,neumann,nielsen,renes}.

We recommend clasical monographs  \cite{bhatia,hardy,marshall} as collections of the most important  inequalities and results in the finite-dimensional majorization theory with their applications onto various fields in mathematics.  See also \cite{ando,hasani,mirsky,von}. 

Notations, preliminaries and some  published results which are used in this paper are contained in Section 2. Section 3   provides some useful properties of introduced increasable doubly stochastic operators and  submajorization relation on $\ell^p(I)^+$, defined by  $f\prec_s g$  whenever $f=Dg$  for some increasable doubly substochastic operator $D$. This relation represents an extension of  the  submajorization between two $n$ dimensional vectors \cite[Theorem I.2.C.4]{marshall} which contains positive real numbers. We provide that, in some sense, the submajorization relation may be considered as a partial order on the positive cone $\ell^p(I)^+$, for some $p\in[1,\infty)$.

In   Section 4 we analyze bounded linear operators on $\ell^p(I)$ which preserve submajorization relation $\prec_s$. Actually, we split the problem on two cases, the submajorization on $\ell^p(I)^+$ when $p\in(1,\infty)$ and the submajorization on $\ell^1(I)^+$. Some complex and technically demanding lemmas and theorems are proved in order to provide and present  concrete forms of linear preservers of submajorization $\prec_s$ (Corollaries \ref{shape} and \ref{shape2}).

\section{Notations and Preliminaries}

In this section we will present notations and the most important results which we will use in the paper.

Let $f:I\longrightarrow
\mathbb{R}$ be an arbitrary function, where $I$ is a non-empty set. The function $f$ is summable if there exists a real number
$\sigma$ with the following property:\\
For every $\epsilon > 0$  there is a finite set $J_0\subseteq I$ such
that
$$\left| \sigma -\sum_{j\in J}f(j)\right|\leq \epsilon$$
whenever $J\subseteq I$ is a finite set and $J_0\subseteq J$. In this  case, such a $\sigma$  is unique and
called the sum of $f$ and we use denotation $\sigma=\sum_{i\in
I}f(i)$.

The Banach space of all functions $f:I\longrightarrow
\mathbb{R}$, where $I$ is non-empty set  and $p\in
[1,\infty)$  with the property $\sum_{i\in I}|f(i)|^p< \infty$ is considered,  and  it  is denoted by $\ell^p(I)$. This space is equipped with $p$-norm 
$$\|f\|_p:=\left(\sum_{i\in I}|f(i)|^p\right)^{\frac 1 p}<
\infty.$$ 
Each function $f\in \ell^p(I)$ may be represented  by 
$$f=\sum_{i\in
I}f(i)e_i,$$ 
where  functions $e_i:I\longrightarrow \mathbb{R}$, $i\in I$ are
defined by Kronecker delta, i.e., $e_i(j)=\delta_{ij}$, for all $j\in I$. 
The above representation of function $f$  means  that    for every $\epsilon > 0$  there is a finite set $J_0\subseteq I$ such
that
$$\| f -\sum_{j\in J}f(j)e_j\|_p\leq \epsilon$$
whenever $J\subseteq I$ is a finite set and $J_0\subseteq J$. Since  $f\in\ell^p(I)$, the above is easy to check. 

We will consider  positive elements of the Banach space $\ell^p(I)$, $p\in [1,\infty)$, and denote this positive cone by
$$\ell^p(I)^+:=\left\{f\in \ell^p(I):f(i)\geq 0,\forall i\in I       \right\}.$$
In some situations we will consider $I_f^+$ as a subset of $I$ determined by
$$I_f^+:=\{i\in I:f(i)>0\},$$
for some $f\in \ell^p(I)$.

A $q$ is conjugate (dual) exponent of $p$
if $\frac{1}{p}+\frac{1}{q}=1$, when $p,q\in (1,\infty)$. Moreover, the exponents $1$ and
$\infty$ are considered to be dual exponent to each other.
For all $ g\in \ell^q(I)$, the map $f\longrightarrow \langle
f,g\rangle:=\sum_{i\in I}f(i)g(i)$ defines a bounded linear functional on
$\ell^p(I)$. Since $\ell^p(I)^*$ is isometrically
isomorphic with $\ell^q(I)$, the dual Banach space $\ell^p(I)^*$ can be
identified with $\ell^q(I)$.
The map $\langle\cdot,\cdot\rangle:\ell^p(I)\times\ell
^q(I)\longrightarrow
\mathbb{R}$ defined by $\langle
f,g\rangle=\summ_{i\in I}f(i)g(i)$ is called the dual pairing.  Using the dual pairing
$\langle\cdot,\cdot\rangle$ for functions
 $f\in \ell^p(I)$ and $e_i\in\ell^1(I)\subset \ell^q(I)$, we have the representation
 $$f(i)=\langle f,e_i\rangle, \hspace{0.3cm} \forall i\in I,$$
and
$$f=\sum_{i\in I}\langle f,e_i\rangle e_i.$$

A bounded linear operator $A:\ell^p(I)\rightarrow\ell^p(I)$ may be represented by a matrix $[a_{ij}]_{i,j\in I}$, which may be finite or infinite depends on cardinality of the set $I$.  If we define matrix elements with  $a_{ij}:=\langle Ae_j,e_i\rangle$, $\forall i,j\in I$ we get the matrix representation of operator $A$ in the following way
\begin{eqnarray}\label{matrixrep}
Af(i)=\sum_{j\in I}a_{ij}f(j),\;\;\;\forall i\in I,
\end{eqnarray}
that is,
$$Af=\sum_{i\in I}\left(\sum_{j\in I}a_{ij}f(j)\right)e_i.$$

\begin{theorem}\label{forallp}\cite[Corollary 3.1]{martin5}
Let $\mathbb{A}=\{a_{ij}:i,j\in I\}$ be a family of real numbers. If this family satisfies conditions 
\begin{equation}\label{columnsup}
\sup_{j\in I}\sum_{i\in I}|a_{ij}|<\infty\end{equation} and \begin{equation}\label{rowsup}\sup_{i\in I}\sum_{j\in I}|a_{ij}|<\infty
\end{equation}
 then this family may be considered as a bounded linear operator $A$ on $\ell^p(I)$, for every $p\in [1,\infty]$, defined by \begin{equation}\label{defoperator2}
Af:=\sum_{i\in I}\left(\sum_{j\in I}a_{ij}f(j)\right)e_i.
\end{equation}
\end{theorem}

 Let $p\in [1,\infty)$ and $A:\ell^p(I)\longrightarrow
\ell^p(I)$ be a bounded linear operator, where $I$ is a
non-empty set. The operator $A$ is called
\begin{itemize}
 \item[i)] a {\em permutation}, if there exists a bijection
    $\theta:I\longrightarrow I$ for which $Ae_j=e_{\theta (j)}$, for
    each $j\in I$;
    \item[ii)] {\em a partial permutation} for sets $I_1\subset
     I$ and $I_2\subset I$, if there exists a bijection
    $\theta:I_1\longrightarrow I_2$ for which $Ae_j=e_{\theta (j)}$  whenever $j\in I_1$, and
    $Ae_j=0$, otherwise.
\end{itemize}
The set of all permutations and set of all  partial permutations on $\ell^p(I)$ are  denoted by $P(\ell^p(I))$ and $pP(\ell^p(I))$, respectively.
Row, column and doubly substochastic operators on $\ell^p(I)$   presented in Definition \ref{defsub} are denoted by $RsS(\ell^p(I))$, $CsS(\ell^p(I))$ and $DsS(\ell^p(I))$, respectively. Doubly stochastic  operators form Definition \ref{DsS} we will denote by  $DS(\ell^p(I))$.

{\rm
\begin{definition}\cite[Definition 4.1]{martin2}\label{gsubmaj}
Let $p\in [1,\infty)$. For two functions $f,g\in {\ell^p(I)}^+$, the  function  $f$ is {\em weakly majorized} by $g$, if there exists a
doubly substochastic operator $D\in DsS(\ell^p(I))$  such that
$f=Dg$, and denote it by $f\prec_{w} g$.
\end{definition}
}

\begin{theorem}\label{subantis}\cite[Theorem 4.2]{martin2}
Let $f,g\in \ell^p(I)^+$, $p\in [1,\infty)$. The next statements are
equivalent:
\begin{itemize}
    \item[i)] $f\prec_w g$ and $g\prec_w f$;
    \item[ii)] There exists $P\in pP(\ell^p(I))$ for sets $I_f^+$ and $I_g^+$ such that $g=Pf$.
\end{itemize}
\end{theorem}

{\rm
\begin{definition}\cite[Definition 3.1]{bahrami}
Let $p\in [1,\infty)$. For two functions $f,g\in \ell^p(I)$, the  function  $f$ is {\em  majorized} by $g$, if there exists a
doubly  stochastic operator $D\in DS(\ell^p(I))$  such that
$f=Dg$, and denote it by $f\prec g$.
\end{definition}
}

\begin{theorem}\label{subantis2}\cite[Theorem 4.3]{martin2}
Let $f,g\in \ell^p(I)^+$, $p\in [1,\infty)$. The next statements are
equivalent:
\begin{itemize}
    \item[i)] $f\prec_w g$ and $g\prec f$;
    \item[ii)] There exists $P\in P(\ell^p(I))$  such that $g=Pf$.
\end{itemize}
\end{theorem}

The space $\ell^p(I)$ is an
ordered Banach space under the natural partial ordering on the set
of real valued functions defined on $I$. The operator $A:\ell^p(I)\longrightarrow \ell^p(I)$, $p\in
[1,\infty)$ is
called positive if $Ag\in\ell^p(I)^+$ for every $g\in \ell^p(I)^+$.

An operator $A^*:\ell^q(I)\longrightarrow \ell^q(I)$ is the adjoint
operator of $A:\ell^p(I)\longrightarrow \ell^p(I)$, $p\in
[1,\infty)$, if $\langle Af,g\rangle= \langle f,A^*g\rangle$,
$\forall f\in \ell^p(I)$, $\forall g\in \ell^q(I)$, where $q$ is the conjugate exponent
of $p$.

\begin{definition}\cite{martin3,martin4}\label{defweakpreservers}
Let $p\in[1,\infty)$. A bounded linear operator $T:\ell^p(I)\rightarrow\ell^p(I)$  is called a preserver of weak majorization on $\ell^p(I)^+$, if $T$ preserves the weak majorization relation, that is, $Tf\prec_w Tg$, whenever $f\prec_w g$, where $f,g\in \ell^p(I)^+$. The set of all linear preservers of weak majorization on $\ell^p(I)^+$ is denoted by $\mathcal{P}_{w}(\ell^p(I)^+)$.
\end{definition}

Let  $\theta: I\rightarrow I$ be one-to-one function. Let $P_\theta:\ell^p(I)\rightarrow\ell^p(I)$, $p\in[1,\infty)$ be a bounded linear operator defined in the following way
\begin{eqnarray}\label{defmapsbyonetoone}
P_\theta(f):=\summ_{k\in I}f(k)e_{\theta(k)},\;\;f\in \ell^p(I).
\end{eqnarray}

\begin{theorem}\cite[Theorem 3.5]{martin3}\label{weakpreserver3}
Let $p\in (1,\infty)$, and let $I$ be an infinite set. Suppose that $T:\ell^p(I)\rightarrow\ell^p(I)$ is a bounded linear operator. The following statements are equivalent:
\begin{itemize}
\item[i)]$T\in  \mathcal{P}_{w}(\ell^p(I)^+)$.
\item[ii)]$Te_j\prec_w Te_k$ and $Te_k\prec_w Te_j$, $\forall k,j\in I$, and for each $i\in I$ there is at most one $j\in I$ such that $\langle Te_j, e_i\rangle>0$.
\item[iii)]$T=\summ_{k\in I_0}\lambda_kP_{\theta_k}$, where $(\lambda_k)_{k\in I_0}\in \ell^p(I_0)^+$, $I_0\subset I$ is at most countable, 
\begin{eqnarray}\label{thetaa}
\theta_k\in \Theta
:=\{\theta_k:I\xrightarrow{1-1} I\,|\; k\in I_0,\hspace{0.2cm}\theta_i (I)\cap\theta_j(I)=\emptyset,\hspace{0.2cm} i\neq j\}.
\end{eqnarray}
\end{itemize}
\end{theorem}
We use two classes of operators on $\ell^1(I)$. Let $\mathcal{P}_1(\ell^1(I)^+)$ be the set of all bounded linear operators on $\ell^1(I)^+$ defined by
\begin{eqnarray}\label{T1}
T_1=\summ_{k\in I_0}\lambda_kP_{\theta_k},
\end{eqnarray} 
where $I_0$ is at most a countable subset of $I$, $(\lambda_k)_{k\in I_0}\in \ell^1(I_0)^+$, every $\theta_k$  belongs to a countable family of one-to-one maps with disjont ranges as in \eqref{thetaa}.
We denote by $\mathcal{P}_2(\ell^1(I)^+)$, the set of all bounded linear operators defined by 
\begin{eqnarray}\label{T2}
T_h(f):=h\sum_{i\in I}f(i),\;\;\forall f\in\ell^1(I),
\end{eqnarray}
where $h\in \ell^1(I)^+$. It is easy to see that each "row" of the operator $T_h$, for example $k\in I$, contains the same elements $h(k)$. In the other words, all "columns" of the operator $T_h$ are the same and they are equal with function $h$.
These two classes of operators are introduced in papers \cite{bahrami,martin4}  and operators in these classes preserve the  weak majorization on $\ell^1(I)^+$ which is provided in the next result.
\begin{theorem}\cite[Theorem 3.3]{martin4}\label{1weakpreserver}
Let $T:\ell^1(I)\rightarrow\ell^1(I)$ be a bounded linear operator, where $I$ is an infinite set. The following statements are equivalent:
\begin{itemize}
\item[i)] $T\in  \mathcal{P}_{w}(\ell^1(I)^+)$;
\item[ii)] There are operators $T_1\in  \mathcal{P}_1(\ell^1(I)^+)$ and $T_2\in  \mathcal{P}_2(\ell^1(I)^+)$ and disjoint sets $I_1,I_2\subset I$ with $I_1\cup I_2=I$ such that $T=T_1+T_2$ where $T_1,T_2$ are chosen to be
\begin{eqnarray*}
\langle T_1f,e_{i_2}\rangle=\langle T_2f,e_{i_1}\rangle =0,\;\;\forall i_1\in I_1,\;\;\forall i_2\in I_2,\;\;\forall f\in \ell^1(I)^+;
\end{eqnarray*}

\item[iii)] There is an at most a countable set $I_0\subset I$ and there is a family
\begin{eqnarray*}
\Theta
:=\{\theta_k:I\xrightarrow{1-1} I\,|\;k\in I_0,\hspace{0.2cm}\theta_i (I)\cap\theta_j(I)=\emptyset,\hspace{0.2cm} i\neq j\}
\end{eqnarray*}
of one-to-one maps, $\theta_k\in \Theta$, $\forall k\in I_0$, and $(\lambda_i)_{i\in I_0}\in \ell^1(I_0)^+$ such that
$$T=\summ_{k\in I_0}\lambda_k P_{\theta_k}+T_h,$$
where $T_h(f):=h\sum_{k\in I}f(k)$, for $h\in \ell^1(I)^+$ with $\langle h,e_j\rangle =0$, $\forall j\in \bigcup_{i\in I_0}\theta_i(I)$;
\item[iv)] $Te_j\prec_w Te_k$ and $Te_k\prec_w Te_j$, $\forall k,j\in I$, and for each $i\in I$, either there exists exactly one $j\in I$ with $\langle Te_j,e_i\rangle>0$ or the set $\{\langle Te_j,e_i\rangle|j\in I\}$ is a singleton.
\end{itemize}
\end{theorem}

\begin{lemma}\cite[Lemma 3.1]{martin4}\label{lemmathereis}
Let $u=\{u_j\}\in\mathbb{R}^n$ and let $\{u_{ij}\,|\,i\in I_0,\,j=1,\ldots ,n\} $ be a family of real numbers, where $I_0$ is at most a countable set. If
\begin{eqnarray}\label{eqlemmathereis}
\sum_{j=1}^n\alpha_ju_j\in\left\{\sum_{j=1}^n\alpha_ju_{ij}\,|\,i\in I_0\right\},
\end{eqnarray}
for all $\alpha=(\alpha_1, \alpha_2,\ldots ,\alpha_n)$ with $\alpha_j>0$ for each $ j=1,\ldots ,n$, then there exists $k\in I_0$ such that $u_j=u_{kj}$, for each $j=1,\ldots,n$.
\end{lemma}

\section{Increasable doubly substochastic operators and submajorization}

At the beginning of this section,   we will introduce the notion of increasable doubly substochastic operators,  based on the relation \eqref{doublysubstochasticmatrix} for the finite-dimensional case.

\begin{definition}\label{defsubbystochastic}
 Let $p\in [1,\infty)$ and let $A:\ell^p(I)\longrightarrow
\ell^p(I)$ be a positive bounded linear operator, where $I$ is a
non-empty set. The operator $A$ is called
 {\em increasable doubly substochastic}, if there is doubly stochastic operator $A_1:\ell^p(I)\longrightarrow
\ell^p(I)$ such that 
\begin{eqnarray}\label{increasable}
\forall i\in I,\;\;\forall j\in I,\;\;\;\langle Ae_j,e_i\rangle\leq \langle A_1e_j,e_i\rangle.\end{eqnarray} 
\end{definition}

\begin{remark}{\rm
In Definition \ref{defsub} we   notice that doubly substochastic operators are defined  from $\ell^p(J)$ to $ \ell^p(I)$, where $I\not = J$ in general. Let $A$   be a bounded linear operator form  $\ell^p(J)$ to $ \ell^p(I)$ which satisfies \eqref{increasable}, where $A_1:\ell^p(J)\longrightarrow
\ell^p(I)$ is a doubly stochastic. 
We recall that for    doubly stochastic operator  $A_1$  have to be $\text{card}(I)=\text{card}(J)$, by \cite[Theorem 2.2]{bahrami}, so  we may choose $I=J$. Because of this,  increasable doubly substochastic operators in Definition \ref{defsubbystochastic} are defined   from $\ell^p(I)$ to $\ell^p(I)$.}
\end{remark}

Increasable doubly substochastic operators on $\ell^p(I)$ introduced in the above Definition \ref{defsubbystochastic}  will be denoted by  $iDsS(\ell^p(I))$.
It is easy to see that 
$$iDsS(\ell^p(I))\subseteq DsS(\ell^p(I))=RsS(\ell^p(I))\cap CsS(\ell^p(I)). $$
When $I$ is a finite set, the equality  $iDsS(\ell^p(I))= DsS(\ell^p(I))$ holds by Theorem \ref{vonnn}. When $I$ is infinite, 
$iDsS(\ell^p(I))\subsetneq DsS(\ell^p(I))$ holds.  Left and right shift operators (presented in \eqref{LR}) are good examples which provide that $iDsS(\ell^p(\mathbb{N}))$ is a proper subset of $DsS(\ell^p(\mathbb{N}))$.

Furthermore,  the norm of doubly substochastic operators is less then or equal to $1$, by \cite[Lemma 3.3]{martin2}. 
Using Definition \ref{defsubbystochastic}, the next lemma is straightforward. 
\begin{lemma}\label{decomp}
Let $p\in[1,\infty)$. For each increasable doubly substochastic operator $D\in iDsS(\ell^p(I))$ there are two operators  $D_1\in DS(\ell^p(I))$ and $D_2\in DsS(\ell^p(I))$ such that 
\begin{eqnarray*}
D_1=D+D_2.
\end{eqnarray*}
\end{lemma}

\begin{lemma}\label{composition} 
Let $p\in[1,\infty)$. The set $iDsS(\ell^p(I))$ is closed under the composition.
\end{lemma}
\begin{proof}

 Let $A,B\in iDsS(\ell^p(I))$. Clearly, $AB$ is a positive operator.
Using Definition \ref{defsubbystochastic} there are corresponding operators $A_1,B_1\in DS(\ell^p(I))$ such that $$\forall i\in I,\;\;\forall j\in I,\;\;\;\langle Ae_j,e_i\rangle\leq \langle A_1e_j,e_i\rangle\;\;\text{and}\;\;\langle Be_j,e_i\rangle\leq \langle B_1e_j,e_i\rangle.$$ 
Now, for arbitrary chosen $i,k\in I$ we get
\begin{eqnarray}\label{closedd}
\langle ABe_k,e_i\rangle&=&\langle A(Be_k),e_i\rangle=A(Be_k)(i)=\sum_{j\in I}\langle Ae_j,e_i\rangle\langle Be_k,e_j\rangle\nonumber\\&\leq& \sum_{j\in I}\langle A_1e_j,e_i\rangle\langle B_1e_k,e_j\rangle=A_1(B_1e_k)(i)\nonumber\\
&=&\langle A_1B_1e_k,e_i\rangle.
\end{eqnarray}
Since $DS(\ell^p(I))$ is closed under the composition by \cite[Theorem 2.4]{bahrami} we obtain $A_1B_1\in DS(\ell^p(I))$,  so  $AB\in iDsS(\ell^p(I))$, by \eqref{closedd}.

\end{proof}

\begin{theorem}
Let $p\in[1,\infty)$. The set $iDsS(\ell^p(I))$ is convex.
\end{theorem}
\begin{proof}
Let $A,B\in iDsS(\ell^p(I))$ and suppose that $A_1,B_1\in DS(\ell^p(I))$ such that $\langle Ae_j,e_i\rangle\leq \langle A_1e_j,e_i\rangle$ and $\langle Be_j,e_i\rangle\leq \langle B_1e_j,e_i\rangle$, for all $ i,j\in I$. 
Let $C=tA+(1-t)B$. Clearly, $C$ is a positive operator. 
Furthermore,
\begin{eqnarray}
\langle Ce_j,e_i\rangle&=&t\langle Ae_j,e_i\rangle+(1-t)\langle Be_j,e_i\rangle\nonumber\\&\leq& t\langle A_1e_j,e_i\rangle+(1-t)\langle B_1e_j,e_i\rangle=\langle C_1e_j,e_i\rangle,\,\;\forall i,j\in I,\nonumber
\end{eqnarray}
where $C_1=tA_1+(1-t)B_1$. It is easy to see that $C_1\in DS(\ell^p(I))$  because $DS(\ell^p(I))$ is a convex set, by \cite[Theorem 3.3]{martin}.
Thus, $C\in iDsS(\ell^p(I))$, so it follows that $iDsS(\ell^p(I))$ is a convex set. 


\end{proof}

\begin{definition}\label{gsubmaj2}
Let $p\in [1,\infty)$. For two functions $f,g\in {\ell^p(I)}^+$, the function $f$ is {\em  submajorized} by $g$, if there exists a increasable
doubly substochastic operator $D\in iDsS(\ell^p(I))$  such that
$f=Dg$, which is denoted by $f\prec_{s} g$.
\end{definition}

Clearly, $f\prec_{s} g$ implies $f\prec_{w} g$, by $iDsS(\ell^p(I))\subset DsS(\ell^p(I))$. The opposite direction is not true in general.
 
It is easy to check using Lemma \ref{decomp} that for two fixed functions $f,g\in\ell^p(I)^+$ the relation $f\prec_s g $ holds if and only if there exists $h\in \ell^p(I)^+$ such that $f\leq h$ and $h\prec g$, where $\leq $ is the  entrywise order ($f(i)\leq h(i),\;\;\forall i\in I$).
 
\begin{example}{\rm Let $p\in[1,\infty)$, $g=(g_1,g_2,g_3,\ldots ,g_n,\ldots)\in\ell^p(\mathbb{N})$, where $g_i>0$, $\forall i\in \mathbb{N}$ and $g_i<g_j$, whenever $i>j$. Let 
$$f:=Rg=(0,g_1,g_2,g_3,\ldots ,g_n,\ldots),$$ where $R$ is the right shift operator defined in \eqref{LR}. Clearly $f\prec_w g$, since $R\in DsS(\ell^p(I))$.
We claim that $f\not\prec_{s}g$. Suppose that $f=Dg$ holds   for an arbitrary chosen doubly substochastic operator $D:\ell^p(\mathbb{N})\rightarrow\ell^p(\mathbb{N})$. Now,
$$0=f_1=\sum_{j=1}^\infty \langle De_j,e_1\rangle g_j$$
so we get $\langle De_j,e_1\rangle=0$, $\forall j\in\mathbb{N}$,  because $g$ is a strictly decreasing sequence with non-zero elements. Further, 
$$g_1=f_2=\sum_{j=1}^\infty \langle De_j,e_2\rangle g_j\leq\sum_{j=1}^\infty \langle De_j,e_2\rangle g_1\leq  g_1,$$ 
so the last inequality holds only when $\langle De_1,e_2\rangle=1$ and 
$\langle De_j,e_2\rangle=0$ whenever $j\not=1.$ Because $D\in DsS(\ell^p(I))\subset CsS(\ell^p(I))$ it follows that $\langle De_1,e_i\rangle=0$, whenever $i\not=2$. Now,
$$g_2=f_3=\sum_{j=1}^\infty \langle De_j,e_3\rangle g_j=0+\sum_{j=2}^\infty \langle De_j,e_3\rangle g_j\leq\sum_{j=2}^\infty \langle De_j,e_3\rangle g_2\leq  g_2.$$ Similarly as above, we get that $\langle De_2,e_3\rangle=1$ and 
$\langle De_j,e_3\rangle=0$ whenever $j\not=1$, and since $D\in CsS(\ell^p(I))$ it follows that $\langle De_2,e_i\rangle=0$, whenever $i\not=3$. Continuing this process we obtain that
have to be $D=R\not\in iDsS(\ell^p(I))$. Thus $f\not\prec_{s} g$.
}
\end{example}

\begin{lemma}\label{lemasubantis3}
Let $p\in [1,\infty)$ and let $f,g\in \ell^p(I)^+$. If
 $f\prec_{s} g$ and $g\prec_{w} f$, then
$f\prec g$.
\end{lemma}

\begin{proof}

Let $\{I_f^n:n\in
\mathbb{N}\}$ be a family of disjoint finite subsets of $I_f^+$ for arbitrary chosen $f\in \ell^p(I)^+$, $p\in[1,\infty)$ defined by
$$I_f^1:=\left\{i\in I_f^+:
f(i)=max\{f(j):j\in I_f^+\}\right\}$$ and
$$I_f^n:=\left\{i\in I_f^+:
f(i)=max\left\{f(j):j\in I_f^+\backslash\bigcup_{k=1}^{n-1}I_f^k\right\}\right\}$$
whenever $n\geq 2$. Above maximums exist by the definition of $\ell^p(I)^+$. Obviously, $I_f^+=\bigcup_{k=1}^\infty I_f^k$. If
$I_ f^k\neq\emptyset$, for some $k\in \mathbb{N}$ then we define
$f_k:=f(j)$, for some $j\in I^k_f$. Otherwise $f_k:=0$.

Let  $f\prec_{s} g$ and $g\prec_w f$. Since $f\prec_{s} g$ implies  $f\prec_{w} g$,  there exists $P\in pP(\ell^p(I))$ for sets $I_f^+$ and $I_g^+$ such that $f=Pg$, by Theorem \ref{subantis}. 
Because of this, it is easy to see that 
\begin{equation}\label{equal}
    f_i=g_i \hspace{0.3cm}and\hspace{0.5cm} \text{card}(I^i_f)=\text{card}(I^i_g),
    \hspace{0.5cm}\forall i\in \mathbb{N}.
\end{equation}

Since $f\prec_{s} g$,  there is $D\in iDsS(\ell^p(I))$  such that $f=Dg$ so using Lemma \ref{decomp} there exist two operators  $D_1\in DS(\ell^p(I))$ and $D_2\in DsS(\ell^p(I))$ such that  $D_1=D+D_2$. We claim that 
\begin{eqnarray}\label{D_2g=0}
D_2g=0.
\end{eqnarray}

Let $i\in I_f^1$. Then
\begin{eqnarray*}
  f_1 &=& f(i)=\sum_{j\in I}g(j) De_j(i)\\
  &=&\sum_{j\in I_g^1}g_1
De_j(i) +\sum_{j\in I\setminus I_g^1}g(j) De_j(i) \\
 &\leq& \sum_{j\in
I_g^1}g_1De_j(i) +\sum_{j\in I\setminus
I_g^1}g_1 De_j(i)\leq g_1=f_1.
\end{eqnarray*}
It follows that $\sum_{j\in I_g^1} 
De_j(i)=1$ and $\sum_{j\in I\setminus
I_g^1} De_j(i)=0$, by $D\in iDsS(\ell^p(I))$.  Thus, for every $i\in I_f^1$ we obtain $$D_1e_j(i)=De_j(i)\;\;\text{and}\;\;\; D_2e_j(i)=0,\;\;\;\forall j\in I.$$  

Since,
$$\text{card}(I_g^1)=\text{card}(I_f^1)=\sum_{i\in I_f^1}\sum_{j\in I_g^1} De_j(i)=\sum_{j\in I_g^1}\sum_{i\in I_f^1} De_j(i),$$
we get $\sum_{i\in I_f^1} De_j(i)=1$, $\forall j\in I_g^1$ and $\sum_{i\in I\setminus I_f^1} De_j(i)=0$, $\forall j\in I_g^1$. Thus, $\forall j\in I_g^1$ we have
$$D_1e_j(i)=De_j(i) \;\;\text{and}\;\;\;D_2e_j(i)=0,\;\;\;\forall i\in I.$$

Let $k\in I_f^{2}$. Using above facts we obtain
\begin{eqnarray*}
  f_{2} &=& f(k)=\sum_{j\in I}g(j) De_j(k)\nonumber\\
  &=&\sum_{j\in
  I_g^{2}}g_{2}
De_j(k) +\sum_{j\in  I_g^1 }g(j) De_j(k)+\sum_{j\in I\setminus  \{I_g^1\cup I_g^2\} }g(j) De_j(k) \\
  &=&\sum_{j\in
  I_g^{2}}g_{2}
De_j(k) +0+\sum_{j\in I\setminus  \{I_g^1\cup I_g^2\} }g(j) De_j(k) \\
  &\leq&\sum_{j\in
  I_g^{2}}g_{2}
De_j(k) +\sum_{j\in I\setminus  \{I_g^1\cup I_g^2\} }g_2 De_j(k) \leq g_{2}=f_2.
\end{eqnarray*}
It follows that $\sum_{j\in I_g^2} 
De_j(k)=1$ and $\sum_{j\in I\setminus I_g^2} 
De_j(k)=0$, by $D\in iDsS(\ell^p(I))$.    
Thus,
$$D_1e_j(k)=De_j(k)\;\;\text{and}\;\;\; D_2e_j(k)=0,\;\;\;\forall j\in I,$$  
for each $k\in I_f^{2}$.

Similarly as above, using $\text{card}(I_g^2)=\text{card}(I_f^2)$ and changing the order of summation, we get
$$\text{card}(I_g^2)=\sum_{j\in I_g^2}\sum_{i\in I_f^2} De_j(i),$$
so $\sum_{i\in I_f^2} De_j(i)=1$, $\forall j\in I_g^2$ and $\sum_{i\in I\setminus I_f^2} De_j(i)=0$, $\forall j\in I_g^2$. Thus $\forall j\in I_g^2$ we have
$$D_1e_j(i)=De_j(i) \;\;\text{and}\;\;\;D_2e_j(i)=0,\;\;\;\forall i\in I.$$ 

Continuing this process, we get  for arbitrary chosen $n\in\mathbb{N}$  that
\begin{eqnarray}\label{l1}
D_2e_j(k)=0,\;\;\;\forall k\in I_f^n,\;\forall j\in I.
\end{eqnarray}
and
\begin{eqnarray}\label{l2}
D_2e_j(i)=0,\;\;\;\forall j\in I_g^n,\;\forall i\in I, 
\end{eqnarray}
hold.
Finally, using \eqref{l2} we obtain
$$D_2g(i)=\sum_{j\in I}g(j)D_2e_j(i)=\sum_{j\in I^+_g}g(j)D_2e_j(i)=\sum_{n=1}^{\infty}\sum_{j\in I^n_g}g(j)D_2e_j(i)=0,$$ for each $i\in I$.
Now,
$f=Dg=(D_1-D_2)g=D_1g$, that is $f\prec g$.

\end{proof}

As direct consequence of Lemma \ref{lemasubantis3}, we get the following three corollaries.
\begin{corollary}\label{subantis3}
Let $f,g\in \ell^p(I)^+$, $p\in [1,\infty)$. The next statements are
equivalent:
\begin{itemize}
    \item[i)] $f\prec_{s} g$ and $g\prec_{w} f$;
    \item[ii)] There exists $P\in P(\ell^p(I))$ such that $g=Pf$.
\end{itemize}
\end{corollary}
\begin{proof}
Let $f\prec_{s} g$ and $g\prec_{w} f$. Using Lemma \ref{lemasubantis3} we get that $f\prec g$. The rest follows by Theorem \ref{subantis2}.

Conversely, if $g=Pf$ for $P\in P(\ell^p(I))\subset iDsS(\ell^p(I))$ it follows that $f=P^{-1}g$ where $P^{-1}\in P(\ell^p(I))\subset DsS(\ell^p(I))$, so we get $f\prec_{s} g$ and $g\prec_{w} f$.
\end{proof}

The next example shows that condition $f\prec_{s} g$ in the above corollary can not be replaced by $f\prec_{w} g$.
\begin{example}{\rm
Let $f\in \ell^1(\mathbb{N})$ defined by $f(i)=1/i^2$, $i\in\mathbb{N}$ that is, $f=(1,\frac{1}{4},\frac{1}{9},\frac{1}{16}\ldots )$. Using right shift operator $R$, we get that $g:=Rf=(0,1,\frac{1}{4},\frac{1}{9},\frac{1}{16}\ldots)$. Similarly, $f=Lg$. It follows that $f$ and $g$ are mutually weakly majorized so they are partial permutations of each other by Theorem \ref{subantis}. However, $g(1)=0\not\in f(I)$. Thus there is no permutation $P\in\ell^1(I)$ to be $f=Pg$.

}\end{example}

Since $g\prec_s f$ implies $g\prec_w f$, using Corollary \ref{subantis3} we obtain the following  result.
\begin{corollary}\label{subantis4}
Let $f,g\in \ell^p(I)^+$, $p\in [1,\infty)$. The next statements are
equivalent:
\begin{itemize}
    \item[i)] $f\prec_{s} g$ and $g\prec_{s} f$;
    \item[ii)] There exists $P\in P(\ell^p(I))$ such that $g=Pf$.
\end{itemize}
\end{corollary}
As corollary of the above result we obtain an analogue of \cite[Theorem 3.5]{bahrami} for positive cone $\ell^p(I)^+$.
\begin{corollary}\label{subantis5}
Let $f,g\in \ell^p(I)^+$, $p\in [1,\infty)$. The next statements are
equivalent:
\begin{itemize}
    \item[i)] $f\prec g$ and $g\prec f$;
    \item[ii)] There exists $P\in P(\ell^p(I))$ such that $g=Pf$.
\end{itemize}
\end{corollary}

\begin{corollary}
The submajorization relation "$\prec_{s}$"
when $p\in[1,\infty)$, is reflexive and transitive relation i.e.
"$\prec_{s}$" is a pre-order. If we identify all
functions which are permutations of each other, then we may
consider "$\prec_{s}$" as a partial order.
\end{corollary}
\begin{proof}
Reflexivity is straightforward. Transitivity follows from Lemma \ref{composition}. If we identify all
functions which are  permutations of each other, then relation $\prec_{s}$ is antisymmetric, by Corollary \ref{subantis4}.
\end{proof}



\section{Linear preservers of submajorization when $I$ is an infinite set}

We reformulate the notion of linear preservers of submajorization relation on $\ell^p(I)$.

\begin{definition}\label{defweakpreservers}
A bounded linear operator $T:\ell^p(I)\rightarrow\ell^p(I)$ is called a preserver of  submajorization on $\ell^p(I)^+$, if $T$ preserves the submajorization relation, that is, $Tf\prec_{s} Tg$, whenever $f\prec_{s} g$, where $f,g\in \ell^p(I)^+$. The set of all linear preservers of  submajorization on $\ell^p(I)^+$ is denoted by $\mathcal{P}_{s}(\ell^p(I)^+)$.
\end{definition}

\begin{theorem}\label{existsubstoch}
Let $D\in iDsS(\ell^p(I))$, $p\in [1,\infty)$, and suppose that
\begin{eqnarray}\label{deftheta}
\Theta :=\{\theta_k:I\xrightarrow{1-1} I\,|\;k\in I_0,\hspace{0.2cm}\theta_i (I)\cap\theta_j(I)=\emptyset,\hspace{0.2cm} i\neq j\}
\end{eqnarray}
is a family of one-to-one maps on $I$ with disjoint images, where $I_0$ is at most a countable set. Then there is at least one $S\in iDsS(\ell^p(I))$ such that $P_\theta D=SP_\theta$, $\forall \theta\in \Theta$.
\end{theorem}

\begin{proof}
Let $D\in iDsS(\ell^p(I))$. There are two operators $D_1\in DS(\ell^p(I))$ and $D_2\in DsS(\ell^p(I))$ such that 
\begin{eqnarray}\label{decomp2}
D_1=D+D_2
\end{eqnarray}
by Lemma \ref{decomp}.
Now, using \cite[Lemma 4.2]{bahrami}  there exists operator $S_1\in DS(\ell^p(I))$ such that $P_\theta D_1=S_1 P_\theta$, $\forall \theta \in \Theta$. Actually, we can see in the proof of above mentioned theorem that operator $S_1$ is defined by
\begin{equation}
\langle S_1e_j,e_i\rangle=\left\{\begin{matrix}
\langle D_1e_{\theta^{-1}(j)},e_{\theta^{-1}(i)}\rangle, & \;\;\; i,j\in \theta(I)\;\;\text{for some}\;\;\theta\in \Theta, \\
\vspace{0.1cm}
1,  & \hspace{0.2cm}i,j\not\in \cup_{\theta\in\Theta}\theta(I)\;\;\text{and}\;\; j= i,\\
\vspace{0.1cm}
0,  &   \text{otherwise.}\\
\end{matrix}
\right.\nonumber
\end{equation}

In the similar way, using \cite[Theorem 3.2]{martin3} there is an operator $S_2\in DsS(\ell^p(I))$ such that $P_\theta D_2=S_2 P_\theta$, $\forall \theta \in \Theta$, defined by
\begin{equation}
\langle S_2e_j,e_i\rangle=\left\{\begin{matrix}
\langle D_2e_{\theta^{-1}(j)},e_{\theta^{-1}(i)}\rangle, & \;\;\;  i,j\in \theta(I),\;\;\text{for some}\;\;\theta\in \Theta, \\
\vspace{0.1cm}
a,  & \hspace{0.2cm}i,j\not\in \cup_{\theta\in\Theta}\theta(I), \;\;\text{and}\;\;  j= i,\\
\vspace{0.1cm}
0,  & \text{otherwise,}\\
\end{matrix}
\right.\nonumber
\end{equation}
where $0\leq a\leq 1$.
 Clearly, operator $S_2$ is not uniquely determined. 
We define bounded linear operator $S:=S_1-S_2$. The operator $S$ has form
\begin{eqnarray}\label{operatorS}
\langle Se_j,e_i\rangle=\left\{\begin{matrix}
\langle D_1e_{\theta^{-1}(j)},e_{\theta^{-1}(i)}\rangle-\langle D_2e_{\theta^{-1}(j)},e_{\theta^{-1}(i)}\rangle, & \; i,j\in \theta(I) \;\text{for some}\;\theta\in \Theta, \\
\vspace{0.1cm}
1-a,  & \hspace{0.2cm}i,j\not\in \cup_{\theta\in\Theta}\theta(I)  \;\;\text{and}\;\;  j= i,\\
\vspace{0.1cm}
0,  & \text{otherwise.}\\
\end{matrix}
\right.
\end{eqnarray}
Obviously, $S\in iDsS(\ell^p(I))$ by the above representation \eqref{operatorS} and the decomposition \eqref{decomp2}. Now,
$$P_\theta D=P_\theta(D_1-D_2)=P_\theta D_1-P_\theta D_2=S_1P_\theta-S_2P_\theta=SP_\theta.$$

\end{proof}

\begin{theorem}\label{subsetweak}
Let $I$ be an infinite set and let $p\in (1,\infty)$. Then, $\mathcal{P}_{w}(\ell^p(I)^+)\subset \mathcal{P}_{s}(\ell^p(I)^+)$ holds.
\end{theorem}
\begin{proof}
Let $T\in\mathcal{P}_{w}(\ell^p(I)^+)$. Using Theorem \ref{weakpreserver3} we get
$$T=\summ_{k\in I_0}\lambda_kP_{\theta_k},$$ where $(\lambda_k)_{k\in I_0}\in \ell^p(I_0)^+$, $I_0\subset I$ is at most countable and for each $k\in I_0$, $\theta_k\in \Theta
=\{\theta_k:I\xrightarrow{1-1} I\,|\; k\in I_0,\hspace{0.2cm}\theta_i (I)\cap\theta_j(I)=\emptyset,\hspace{0.2cm} i\neq j\}$.

Let $f\prec_{s} g$. There is an operator $D\in iDsS(\ell^p(I))$ such that $f=Dg$. Using Theorem \ref{existsubstoch}, there is an operator $S\in iDsS(\ell^p(I))$ such that $P_\theta D=SP_\theta$, $\forall \theta\in \Theta$. Hence,
\begin{eqnarray}
Tf&=&\summ_{k\in I_0}\lambda_kP_{\theta_k}(f)=\summ_{k\in I_0}\lambda_kP_{\theta_k}(Dg)\nonumber\\
&=&\summ_{k\in I_0}\lambda_kSP_{\theta_k}(g)=S\left(\summ_{k\in I_0}\lambda_kP_{\theta_k}(g)\right)=S(Tg)\nonumber,
\end{eqnarray}
thus $Tf\prec_{s} Tg$, so $T\in \mathcal{P}_{s}(\ell^p(I)^+)$.
\end{proof}

In the sequel, we will show that every linear preserver of the  submajorization $(\prec_{s})$ preserve the weak majorization $(\prec_{w})$, when $p\in(1,\infty)$ and when $I$ is an infinite set.

\begin{theorem}\label{jednauvrsti}
Let  $T\in \mathcal{P}_{s}(\ell^p(I))$ where $I$ is an infinite set and $p\in (1,\infty)$. Then for each pair of distinct elements $j_1,j_2\in I$, functions $Te_{j_1}$ and $Te_{j_2}$ are permutations of each other and $Te_{j_1}(i)\cdot Te_{j_2}(i)=0$, for all $i\in I$.
\end{theorem}

\begin{proof}
Since $e_{j_1}\prec_{s}e_{j_2}$ and $e_{j_2}\prec_{s}e_{j_1}$, we have $Te_{j_1}\prec_{s}Te_{j_2}$ and $Te_{j_2}\prec_{s}Te_{j_1}$. Using Corollary \ref{subantis4}, functions  $Te_{j_1}$ and $Te_{j_2}$ are permutations of each other.

In order to show second part,  we will suppose contrary that there are two elements $j_1,j_2\in I$ and there exists $k\in I$ such that $Te_{j_1}(k)\cdot Te_{j_2}(k)\not =0$. Let $c_1=Te_{j_1}(k) \not=0$ and $c_2=Te_{j_2}(k)\not=0$.
Since the "column" $Te_{j_1}$ in the matrix form of operator $T$  is in $\ell^p(I)$, we have $\lim_{k\rightarrow\infty}Te_{j_1}(k)=0$. Hence, there is a finite set   defined by
 $$\mathcal{C}_1:=\{i\in I:Te_{j_1}(i)=c_1\}.$$
 Clearly, $\mathcal{C}_1$ is a  non-empty set because $k\in \mathcal{C}_1$.
 
Let $j_3\in I\setminus\{j_1,j_2\}$.
For arbitrary chosen $a>0$ and $b>0$ we have $ae_{j_1}+be_{j_2}\prec_{s} ae_{j_1}+be_{j_3}$ and $ae_{j_1}+be_{j_3}\prec_{s} ae_{j_1}+be_{j_2}$. Therefore, $aTe_{j_1}+bTe_{j_2}\prec_{s} aTe_{j_1}+bTe_{j_3}$ and $aTe_{j_1}+bTe_{j_3}\prec_{s} aTe_{j_1}+bTe_{j_2}$. Now, using Corollary \ref{subantis4}, functions $aTe_{j_1}+bTe_{j_2}$ and $aTe_{j_1}+bTe_{j_3}$ are permutations of each other, that is
$$ac_1+bc_2=aTe_{j_1}(k)+bTe_{j_2}(k)\in\{ aTe_{j_1}(i)+bTe_{j_3}(i)\;:\;i\in I\}.$$
Since the above set is at most countable, using Lemma \ref{lemmathereis} for $n=2$ we obtain that there exists $k_0\in I $ such  that $Te_{j_1}(k_0)=c_1$ and $Te_{j_3}(k_0)=c_2$, so we conclude that $k_0\in \mathcal{C}_1$. Since $j_3$ is arbitrary chosen element of the infinite set $I\setminus\{j_1,j_2\}$ and  since  the set $\mathcal{C}_1$ is finite, there is at least one element $k_1\in \mathcal{C}_1$ for which there is an infinite set 
$\mathcal{T}\subset I$ such that $Te_{j_1}(k_1)=c_1$ and $Te_{t}(k_1)=c_2$, $\forall t\in\mathcal{T}$.
Now, for adjoint operator $T^*:\ell^q(I)\rightarrow\ell^q(I)$ of $T$, where $q$ is conjugate exponent of $p$, we get
\begin{eqnarray*}
\|T^*e_{k_1}\|^q&=&\sum_{j\in I}|T^*e_{k_1}(j)|^q\geq \sum_{t\in \mathcal{T}}|\langle T^*e_{k_1},e_t\rangle|^q=\sum_{t\in \mathcal{T}}|\langle Te_t, e_{k_1}\rangle|^q\\
&=& \sum_{t\in \mathcal{T}}c_2^q=\infty,
\end{eqnarray*}
which is impossible.

In the other words, we conclude that the "row"  $k_1$   of operator $T$ contains infinite nonzero elements  which are mutually equal. It implies that the  same holds for appropriate $k_1$ "column"  of adjoint operator $T^*$. Since $T^*:\ell^q(I)\rightarrow\ell^q(I)$ have to be $T(e_{k_1})\in \ell^q(I)$, which is a contadiction with the above fact.

\end{proof}

Using above two theorems we obtain that preservers of weak majorization $(\prec_{w})$ and submajorization $(\prec_{s})$ on $\ell^p(I)^+$ coincide when  $I$ in an infinite set and when $p\in(1,\infty)$. Using Theorem \ref{weakpreserver3}, we get the following characterization of submajorization preservers on $\ell^p(I)^+$.

\begin{corollary}\label{shape}
Let $p\in (1,\infty)$, and let $I$ be an infinite set. Suppose that $T:\ell^p(I)\rightarrow\ell^p(I)$ is a bounded linear operator. The following statements are equivalent:
\begin{itemize}
\item[i)]$T\in  \mathcal{P}_{s}(\ell^p(I)^+)$.
\item[ii)]$Te_j\prec_{s} Te_k$ and $Te_k\prec_{s} Te_j$, $\forall k,j\in I$, and for each $i\in I$ there is at most one $j\in I$ such that $\langle Te_j, e_i\rangle>0$.
\item[iii)]$T=\summ_{k\in I_0}\lambda_kP_{\theta_k}$, where $(\lambda_k)_{k\in I_0}\in \ell^p(I_0)^+$, $I_0\subset I$ is at most countable, $\theta_k\in \Theta
:=\{\theta_k:I\xrightarrow{1-1} I\,|\; k\in I_0,\hspace{0.2cm}\theta_i (I)\cap\theta_j(I)=\emptyset,\hspace{0.2cm} i\neq j\}.$
\end{itemize}
\end{corollary}

In the sequel, we will find the  proper form of linear preservers of submajorization on $\ell^1(I)$, when $I$ is an infinite set.

\begin{theorem}\label{w_sub_s}
Let $I$ be an infinite set. Then, $\mathcal{P}_{w}(\ell^1(I)^+)\subset \mathcal{P}_{s}(\ell^1(I)^+)$ holds.
\end{theorem}

\begin{proof}
Let $T\in\mathcal{P}_{w}(\ell^1(I)^+)$. Using Theorem \ref{1weakpreserver} statement $ii)$, we get that there is a decomposition of the operator $T$ in the follofing way:
$$T=T_1+T_2,$$
where $T_1\in  \mathcal{P}_1(\ell^1(I)^+)$ and $T_2\in  \mathcal{P}_2(\ell^1(I)^+)$, where sets $I_1,I_2\subset I$ are disjoint with $I_1\cup I_2=I$ and operators $T_1,T_2$ are chosen to be
\begin{eqnarray}\label{uslov}
\langle T_1f,e_{i_2}\rangle=\langle T_2f,e_{i_1}\rangle =0,\;\;\forall i_1\in I_1,\;\;\forall i_2\in I_2,\;\;\forall f\in \ell^1(I)^+.
\end{eqnarray}
Using \eqref{T1}, we have
\begin{eqnarray}
T_1=\summ_{k\in I_0}\lambda_kP_{\theta_k}
\end{eqnarray} 
where $I_0$ is at most a countable subset of $I$, $(\lambda_k)_{k\in I_0}\in \ell^1(I_0)^+$ and  for every $k\in I_0$  we have
$$\theta_k\in \Theta
=\{\theta_k:I\xrightarrow{1-1} I\,|\; k\in I_0,\hspace{0.2cm}\theta_i (I)\cap\theta_j(I)=\emptyset,\hspace{0.2cm} i\neq j\}.$$
If there is $i\in I_0$ such that $\lambda_i=0$ then we will consider the set $I_0\setminus \{i\}$ instead of $I_0$. Because of this, we may assume that $\lambda_j>0$, for every $j\in I_0$.

Suppose that $f\prec_{s} g$ for fixed $f,g\in \ell^1(I)^+$. It follows that there exists $D\in iDsS(\ell^1(I))$ such that $f=Dg$. Using Theorem \ref{existsubstoch} we conclude that there is $S\in iDsS(\ell^1(I))$ such that 
$$P_{\theta_k} D=SP_{\theta_k},\;\;\forall k\in I_0.$$ Obviously, operator $S$ is not unique and it is defined by \eqref{operatorS}, where $0\leq a \leq 1$. Similarly as in Theorem \ref{subsetweak} we obtain that $T_1$ preserve submajorization relation:
\begin{eqnarray}\label{i}
T_1f=T_1Dg=\summ_{k\in I_0}\lambda_kP_{\theta_k}(Dg)
=\summ_{k\in I_0}\lambda_kSP_{\theta_k}(g)=S(T_1g).
\end{eqnarray}

Next, changing the order of summation we obtain 
\begin{eqnarray*}
\|f\|&=&\sum_{i\in I}|f(i)|=\sum_{i\in I}f(i)=\sum_{i\in I}\sum_{j\in I}\langle De_j,e_i\rangle g(j)\nonumber\\
&=&\sum_{j\in I}\sum_{i\in I}\langle De_j,e_i\rangle g(j)=\sum_{j\in I} g(j)\sum_{i\in I}\langle De_j,e_i\rangle\leq \|g\|.
\end{eqnarray*}
Thus, inequality $\|f\|\leq \|g\|$ holds.
If we set $a:=1-\frac{\|f\|}{\|g\|}$, using the above argument we get  that $0\leq a\leq 1.$ Now, using \eqref{T2} we get
\begin{eqnarray}\label{ll}
T_2Dg=T_2f=h\sum_{i\in I}f(i)=h\|f\|=h(1-a)\|g\|=(1-a)T_2g,
\end{eqnarray}
where $h:=T_2e_j$, for some $j\in I$.

We claim that 
\begin{eqnarray}\label{lll}
ST_2=(1-a)T_2.
\end{eqnarray}
 Firstly, we will show that $Se_k=(1-a)e_k$, for every $k\in I_2$. Fix $k\in I_2$. We have that $\langle T_1f,e_k\rangle=0$, by \eqref{uslov}. Since,
\begin{eqnarray*}
\langle T_1f,e_k\rangle&=&\sum_{j\in I}f(j)\langle T_1e_j,e_k\rangle=\sum_{j\in I}f(j)\sum_{i\in I_0}\lambda_i \langle P_{\theta_i}e_j,e_k\rangle\\
&=&\sum_{j\in I}f(j)\sum_{i\in I_0}\lambda_i \langle e_{\theta_i(j)},e_k\rangle
\end{eqnarray*}
It follows that $k\not\in \cup_{i\in I_0}\theta_i(I)$, because $f$ is arbitrary fixed and $\lambda$ is positive. Now, using the definition \eqref{operatorS} of the operator $S$ we get
$Se_k=(1-a)e_k$. Now, using \eqref{uslov} we obtain

\begin{eqnarray*}
ST_2u&=&(\sum_{j\in I}u_j) Sh=(\sum_{j\in I}u_j)\sum_{k\in I_2}h(k)Se_k= (1-a)(\sum_{j\in I}u_j)\sum_{k\in I_2}h(k)e_k\\&=&(1-a)T_2u,
\end{eqnarray*}
for every $u\in \ell^1(I)^+$, so \eqref{lll}  is provided.
Finally, using \eqref{i}, \eqref{ll} and \eqref{lll} we obtain
\begin{equation*}
Tf=(T_1+T_2)Dg=T_1Dg+T_2Dg=ST_1g+(1-a)T_2g=ST_1g+ST_2g=STg.
\end{equation*}
Since, $S\in iDsS(\ell^1(I))$, we get $Tf\prec_{s}Tg$, that is $T\in \mathcal{P}_{s}(\ell^1(I)^+)$.
\end{proof}

In order to show that $\mathcal{P}_{s}(\ell^1(I)^+)$ is a subset of $\mathcal{P}_{w}(\ell^1(I)^+)$ we need the following lemmas.

\begin{lemma}\label{lemmasuperthereis2}
Let $T\in\mathcal{P}_{s}(\ell^1(I)^+)$.  Suppose that $Q$ is a finite subset of $I$ and let  $\Delta:Q\rightarrow Q$ be a bijection.  For every $a\in I$ there is $b\in I$ such that
\begin{eqnarray}\label{eqlemmasuperthereis2}
 \langle Te_i,e_a\rangle= \langle Te_{\Delta(i)},e_b\rangle,\;\;\;\forall i\in Q.
\end{eqnarray}

\end{lemma}

\begin{proof}
 Since,
$$\sum_{i\in Q}a_ie_i\prec_{s} \sum_{i\in Q}a_i e_{\Delta(i)}\;\; {\rm and}\;\; \sum_{i\in Q}a_i e_{\Delta(i)}\prec_{s}  \sum_{i\in Q}a_ie_i$$
we get
$$\sum_{i\in Q}a_iTe_i\prec_{s}  \sum_{i\in Q}a_i Te_{\Delta(i)}\;\; {\rm and}\;\; \sum_{i\in Q}a_i Te_{\Delta(i)}\prec_{s}  \sum_{i\in Q}a_iTe_i$$
for each $(a_{i_1},a_{i_2},\ldots ,a_{i_m})$ with $a_{i_j}>0$ for every $i_j$, where $m=\text{card}(Q)\in \mathbb{N}$.
We get that functions $\summ_{i\in Q}a_iTe_i$ and $\summ_{i\in Q}a_i Te_{\Delta(i)}$ are permutations of each other by Corollary \ref{subantis4}, that is
$$\sum_{i\in Q}a_i \langle Te_i,e_a\rangle\in\left\{\sum_{i\in Q}a_i  \langle Te_{\Delta(i)},e_k\rangle \mid k\in I\right\}.$$
Since codomain of the positive operator $T$ is $\ell^1(I)$ we have $Te_j\in\ell^1(I)^+$ and $\text{card}(\text{Im}(Te_j))\leq \aleph_0$, so $\sum_{i\in Q}a_i    Te_{\Delta(i)}$ is at most a countable set.
We get that for fixed $a\in I$ there is a $b\in I$ such that (\ref{eqlemmasuperthereis2}) holds, by Lemma \ref{lemmathereis}.
\end{proof}

\begin{lemma}\label{lemmaequal}
Let $T\in\mathcal{P}_{s}(\ell^1(I)^+)$, where $I$ is an infinite set. If there are two distinct $n,m\in I$ such that $\langle Te_n,e_r\rangle>0$ and $\langle Te_m,e_r\rangle>0$ for some $r\in I$, then $\langle Te_n,e_r\rangle=\langle Te_m,e_r\rangle$.
\end{lemma}

\begin{proof}
Suppose that there exist $m,n,r\in I$ such that  $\langle Te_m,e_r\rangle>0 $, $\langle Te_n,e_r\rangle>0 $ and $\langle Te_m,e_r\rangle\not =\langle Te_n,e_r\rangle$.

Let $a_1,a_2>0$ and let $l\in I\setminus \{m,n\}$ be arbitrary chosen.
Clearly,
$$a_1e_m+a_2e_n\prec_{s}a_1e_m+a_2e_l\;\;\text{and}\;\;a_1e_m+a_2e_l\prec_{s}a_1e_m+a_2e_n.$$
Because $T\in\mathcal{P}_{s}(\ell^1(I)^+)$ we obtain
$$a_1Te_m+a_2Te_n\prec_{s}a_1Te_m+a_2Te_l\;\;\text{and}\;\;a_1Te_m+a_2Te_l\prec_{s}Ta_1e_m+a_2Te_n.$$
It follows that 
$$a_1 \langle  Te_m,e_r\rangle+a_2 \langle Te_n,e_r\rangle\in\left\{a_1\langle Te_m,e_j\rangle +a_2\langle Te_l,e_j\rangle \; | \; j\in I\right\}.$$
by Corollary \ref{subantis4}.

Since $a_1Te_m+a_2Te_l\in \ell^1(I)^+$, the above set is at most countable, so using Lemma \ref{lemmathereis} for $n=2$ we get that for $l$ there is $k\in I$ such that 
\begin{eqnarray}\label{ttt}
\langle Te_m,e_k\rangle= \langle Te_m,e_r\rangle\;\;\text{and}\;\;\langle Te_l,e_k\rangle= \langle Te_n,e_r\rangle.
\end{eqnarray}

On the other hand, it is clear that for fixed $c\in \mathbb{R}$, $c\not=0$ holds
\begin{eqnarray}\label{aleph0}
\text{card}\{i\in I\;|\;\langle Te_m,e_i\rangle=\text{c}\}<\aleph_0.
\end{eqnarray}
Since $I$ is an infinite set and $l\in I\setminus\{m,n\}$ is arbitrary chosen, using \eqref{ttt} and \eqref{aleph0} there is  $s\in I$ and there is a sequence $ (t_i)_{i\in \mathbb{N}}$ of distinct elements $t_i\in I$ such that 
$$ \langle Te_{m},e_s\rangle=\langle Te_m,e_r\rangle\;\;\text{and}\;\; \langle Te_{t_i},e_s\rangle=\langle Te_n,e_r\rangle,\;\;\forall i\in \mathbb{N}.$$

Let $\{ \Phi _j\} _{j\in\mathbb{N}}$ be a family  where $\Phi_j:=\{t_1,t_2,\ldots t_j\}$ for every $j\in \mathbb{N}$.

 We define  bijections $\gamma_j:\{t_1,t_2,\ldots t_j\}\cup \{m\}\rightarrow \{t_1,t_2,\ldots t_j\}\cup \{m\}$ in the following way:

\begin{equation*}\label{defphisuper}
\gamma_j(x):=\left\{\begin{matrix}
t_j, &  x=m,   \\
m,  & x=t_j  \\
x, &  x\in\{t_1,t_2,\ldots t_{j-1}\}.    \\
\end{matrix}
\right.
\end{equation*}
For each $j\in \mathbb{N}$ there exists $r_j\in I$ such that
\begin{eqnarray}\label{lemmasubequal1}
 \langle Te_{m},e_{r_j}\rangle=\langle Te_{\gamma_j(t_j)},e_{r_j}\rangle =\langle T{e_{t_j}},e_s\rangle= \langle Te_n,e_r\rangle,
\end{eqnarray}
\begin{eqnarray}\label{lemmasubequal2}
 \langle Te_{t_j},e_{r_j}\rangle=\langle Te_{\gamma_j(m)},e_{r_j}\rangle =\langle Te_m,e_s\rangle= \langle Te_m,e_r\rangle,
\end{eqnarray}
by Lemma \ref{lemmasuperthereis2}
Also, for each $x\in\{t_1,t_2,\ldots t_{j-1}\}$ we get
\begin{eqnarray}\label{lemmasubequal3}
\langle Te_{x},e_{r_j}\rangle= \langle Te_{\gamma_j(x)},e_{r_j}\rangle=\langle Te_x,e_s\rangle= \langle Te_{n},e_r\rangle,
\end{eqnarray}
again by Lemma \ref{lemmasuperthereis2}.

If we provide that the set $\{r_j\;|\;j\in \mathbb{N}\}$ is  countable, we obtain
 using \eqref{lemmasubequal1} that  $\langle Te_{m},e_{r_j}\rangle=\langle Te_n,e_r\rangle$, for all $j\in\mathbb{N}$ which   is a contradiction with \eqref{aleph0}.

Let  $k_1<k_2$ for some integers $k_1,k_2 $  and suppose that $ r_{k_1}=r_{k_2}$. Since bijections $\gamma_{k_1}$ and $\gamma_{k_2}$ are different,  using \eqref{lemmasubequal2} for $\gamma_{k_1}$  we obtain  
$$\langle Te_{t_{k_1}},e_{r_{k_1}}\rangle=\langle Te_m,e_r\rangle.$$
Because $k_1<k_2$ we get using  \eqref{lemmasubequal3} for $\gamma_{k_2}$ that 
$$\langle Te_{t_{k_1}},e_{r_{k_1}}\rangle=\langle Te_{t_{k_1}},e_{r_{k_2}}\rangle=\langle Te_{n},e_r\rangle.$$ 
Combine above two facts, we get 
$$\langle Te_m,e_r\rangle=\langle Te_{t_{k_1}},e_{r_{k_1}}\rangle=\langle Te_{n},e_r\rangle$$
 which is a contradiction with the assumption at the beginning of the proof $\langle Te_m,e_r\rangle\not =\langle Te_n,e_r\rangle$. Thus, $ r_{k_1}\not=r_{k_2}$  for all $k_1,k_2\in\mathbb{N}$, so the set $\{r_j\;|\;j\in I\}$ is a countable.

\end{proof}

\begin{theorem}\label{s_sub_w}
Let $I$ be an infinite set. Then, $\mathcal{P}_{s}(\ell^1(I)^+)\subset \mathcal{P}_{w}(\ell^1(I)^+)$ holds.
\end{theorem}

\begin{proof}
Let $T\in\mathcal{P}_{s}(\ell^1(I)^+)$.
Because $e_i\prec_{s} e_j$ and $e_j\prec_{s} e_i$ it follows that $Te_i\prec_{s} Te_j$ and $Te_j\prec_{s} Te_i$. Since, relation $\prec_s$ implies $\prec_w$, we have
$Te_i\prec_{w} Te_j$ and $Te_j\prec_{w} Te_i$, so the first part of     statement  $iv)$ in Theorem \ref{1weakpreserver}  is satisfied.
 
 Suppose that there exist $m,n,r\in I$ such that  $\langle Te_m,e_r\rangle>0 $ and $\langle Te_n,e_r\rangle>0 $.
 Using Lemma \ref{lemmaequal} we get that $\langle Te_m,e_r\rangle=\langle Te_n,e_r\rangle$. Precisely, all non-zero elements in one "row" have to be mutually equal. 
 
We claim that all elements in the "row" indexed by $r$ are the same, that is, there is no zero element. Suppose contrary  that there exists  $l\in I$ such that $\langle Te_l,e_r\rangle=0$. Fix $k\in I\setminus\{m,n,l\}$. We define a bijection 

\begin{equation*}\label{defthetasub}
\theta_k(x):=\left\{\begin{matrix}
m, &  x=m,   \\
n, &  x=n,   \\
k,  & x=l,\\
l, &  x=k.  \\
\end{matrix}
\right.
\end{equation*}

Now, applying Lemma \ref{lemmasuperthereis2} on  bijection $\theta_k$ we get that there exists $i_k\in I$ such that $\langle Te_m,e_{i_k}\rangle=\langle Te_m,e_r\rangle$, $\langle Te_n,e_{i_k}\rangle=\langle Te_n,e_r\rangle$ and $\langle Te_k,e_{i_k}\rangle=\langle Te_l,e_r\rangle=0$.
Using $Te_m,Te_n\in \ell^1(I)^+$, it follows that
\begin{eqnarray}\label{aleph0-2}
\text{card}\{i\in I\;|\;\langle Te_m,e_i\rangle=\langle Te_n,e_i\rangle=\langle Te_m,e_r\rangle\}<\aleph_0.
\end{eqnarray}
It is easy to see that $i_k$ is contained in the above set. Now, since $k$ is arbitrary chosen from infinite set $ I\setminus\{m,n,l\}$,   there exists at least one $s\in I$ ($s$ is contained in the  set considered in \eqref{aleph0-2}) and there is a sequence $ (t_i)_{i\in \mathbb{N}}$ of distinct elements $t_i\in I$ such that $ \langle Te_{m},e_s\rangle=\langle Te_m,e_r\rangle= \langle Te_{n},e_s\rangle=\langle Te_n,e_r\rangle$ and $ \langle Te_{t_i},e_s\rangle=0$, $\forall i\in \mathbb{N}$.
Now, for each $j\in\mathbb{N}$ we define bijections $$\gamma_j:\{t_1,t_2,\ldots t_j\}\cup \{m,n\}\rightarrow \{t_1,t_2,\ldots t_j\}\cup \{m,n\},$$ by
\begin{equation*}\label{defphisuper2}
\gamma_j(x):=\left\{\begin{matrix}
m, &  x=m,   \\
t_j, &  x=n,   \\
n,  & x=t_j,  \\
x, &  x\in\{t_1,t_2,\ldots t_{j-1}\}.    \\
\end{matrix}
\right.
\end{equation*}

Similarly as in \eqref{lemmasubequal1}, \eqref{lemmasubequal2} and \eqref{lemmasubequal3},
for each $j\in I$ and for the suitable bijection $\gamma_j$,  there exists $r_j\in I$ such that
\begin{eqnarray}\label{lemmasubequal4}
 \langle Te_{m},e_{r_j}\rangle = \langle Te_m,e_s\rangle= \langle Te_m,e_r\rangle,
\end{eqnarray}
\begin{eqnarray}\label{lemmasubequal5}
 \langle Te_{t_j},e_{r_j}\rangle=\langle Te_{\gamma_j(n)},e_{r_j}\rangle =\langle Te_n,e_{s}\rangle  =\langle Te_n,e_r\rangle,
\end{eqnarray}
\begin{eqnarray}\label{lemmasubequal6}
 \langle Te_{n},e_{r_j}\rangle=\langle Te_{\gamma_j(t_j)},e_{r_j}\rangle =\langle Te_{t_j},e_s\rangle =0,
\end{eqnarray}
again by Lemma \ref{lemmasuperthereis2}.
Also, for every $x\in\{t_1,t_2,\ldots t_{j-1}\}$ we have
\begin{eqnarray}\label{lemmasubequal7}
\langle Te_{x},e_{r_j}\rangle= \langle Te_{\gamma_j(x)},e_{r_j}\rangle= \langle Te_{x},e_s\rangle= 0.
\end{eqnarray}
by Lemma \ref{lemmasuperthereis2}.

Suppose that $k_1<k_2$. We will show that $r_{k_1}\not=r_{k_2}$. Using \eqref{lemmasubequal7} for $\gamma_{k_2}(x)$ we get that 
\begin{eqnarray*}\label{contrad}
\langle Te_{t_{k_1}},e_{r_{k_2}}\rangle=0.
\end{eqnarray*}
However, using \eqref{lemmasubequal5} for bijection $\gamma_{k_1}(x)$ we obtain 
$$\langle Te_{t_{k_1}},e_{r_{k_1}}\rangle=\langle Te_n,e_r\rangle>0,$$ 
so we get that $\langle Te_{k_1},e_{r_{k_1}}\rangle\not=\langle Te_{k_1},e_{r_{k_2}}\rangle$, therefore $r_{k_1}\not= r_{k_2}$.

Finally, we get that sequence $\{r_j\}_{j\in \mathbb{N}}$ contains mutually district members. Using \eqref{lemmasubequal4} it follows that 
$$\|Te_m\|^p=\sum_{i\in I}\langle Te_m,e_i\rangle\geq\sum_{j\in I}\langle Te_m,e_{r_j}\rangle=\sum_{i=1}^\infty \langle Te_m,e_r\rangle=+\infty,$$ which is impossible.
Thus, there is no  $l\in I$ such that $\langle Te_l,e_r\rangle=0$, so the set $\{\langle Te_j,e_i\rangle\;|\;j\in I\}$ is a singleton.
It follows that operator $T$  satisfies statement  $iv)$ of Theorem \ref{1weakpreserver}, that is $T\in\mathcal{P}_{w}(\ell^1(I)^+)$.

\end{proof}

Using Theorem \ref{w_sub_s} and Theorem \ref{s_sub_w} we obtain that preservers of weak majorization $(\prec_{w})$ and submajorization $(\prec_{s})$ on $\ell^1(I)^+$ coincide when  $I$ in an infinite set. Using Theorem \ref{1weakpreserver}, we get the following characterization of submajorization preservers on $\ell^1(I)^+$.

\begin{corollary}\label{shape2}
Let $T:\ell^1(I)\rightarrow\ell^1(I)$ be a bounded linear operator, where $I$ is an infinite set. The following statements are equivalent:
\begin{itemize}
\item[i)] $T\in  \mathcal{P}_{s}(\ell^1(I)^+)$;
\item[ii)] There are operators $T_1\in  \mathcal{P}_1(\ell^1(I)^+)$ and $T_2\in  \mathcal{P}_2(\ell^1(I)^+)$ and disjoint sets $I_1,I_2\subset I$ with $I_1\cup I_2=I$ such that $T=T_1+T_2$ where $T_1,T_2$ are chosen to be
\begin{eqnarray*}
\langle T_1f,e_{i_2}\rangle=\langle T_2f,e_{i_1}\rangle =0,\;\;\forall i_1\in I_1,\;\;\forall i_2\in I_2,\;\;\forall f\in \ell^1(I)^+;
\end{eqnarray*}

\item[iii)] There is an at most a countable set $I_0\subset I$ and there is a family
\begin{eqnarray*}
\Theta
:=\{\theta_k:I\xrightarrow{1-1} I\,|\;k\in I_0,\hspace{0.2cm}\theta_i (I)\cap\theta_j(I)=\emptyset,\hspace{0.2cm} i\neq j\}
\end{eqnarray*}
of one-to-one maps, $\theta_k\in \Theta$, $\forall k\in I_0$, and $(\lambda_i)_{i\in I_0}\in \ell^1(I_0)^+$ such that
$$T=\summ_{k\in I_0}\lambda_k P_{\theta_k}+T_h,$$
where $T_h(f):=h\sum_{k\in I }f(k)$, for $h\in \ell^1(I)^+$ with $\langle h,e_j\rangle =0$, $\forall j\in \bigcup_{i\in I_0}\theta_i(I)$;
\item[iv)] $Te_j\prec_w Te_k$ and $Te_k\prec_w Te_j$, $\forall k,j\in I$, and for each $i\in I$, either there exists exactly one $j\in I$ with $\langle Te_j,e_i\rangle>0$ or the set $\{\langle Te_j,e_i\rangle|j\in I\}$ is a singleton.
\end{itemize}
\end{corollary}

Every linear preserver $T$ of the submajorization relation ($\prec_s$) on $\ell^p(I)^+$ is a positive operator ($Tf\in\ell^p(I)^+$ whenever $f\in \ell^p(I)^+$). It is an evident consequence of Corollaries \ref{shape} and \ref{shape2}.

We recall that the set of all linear preservers of weak majorization on $\ell^p(I)$, $p\in[1,\infty)$ is a norm-closed by \cite[Theorem 4.4]{martin3} and \cite[Theorem 4.2]{martin4}.
Linear preservers of weak majorization and  submajorization coinside, by Corollaries \ref{shape} and \ref{shape2}, hence the following corollary is straightforward.


\begin{corollary}
Let $I$ be an infinite set, and let $p\in[1,\infty)$. The set $\mathcal{P}_{s}(\ell^p(I)^+)$ is a norm-closed subset of the set of all bounded linear operators on $\ell^p(I)$.
\end{corollary}

\begin{example}\rm{We define maps $\theta_i:\mathbb{N}\rightarrow\mathbb{N}$
 by
\begin{eqnarray}\label{deftet}
\theta_i(j)=i+1+\sum_{k=0}^{i+j-2}k,\;\;\forall i,j\in \mathbb{N}.
\end{eqnarray}

Suppose that there exist $i_1,i_2,j_1,j_2\in\mathbb{N}$ such that
\begin{eqnarray}\label{111}
\theta_{i_1}(j_1)=\theta_{i_2}(j_2).
\end{eqnarray}
It follows that
\begin{equation}\label{100}
i_1+\sum_{k=0}^{i_1+j_1-2}k=i_2+\sum_{k=0}^{i_2+j_2-2}k.
\end{equation}
If $i_1+j_1<i_2+j_2$ then, by \eqref{100}, we get
$$i_1-i_2=\sum_{k=i_1+j_1-1}^{i_2+j_2-2} k\geq i_1+j_1-1.$$
Since $j_1,i_2\in \mathbb{N}$, we get by above that $j_1+i_2\leq 1$ holds, which is impossible. Similarly, if $i_1+j_1>i_2+j_2$ then we get $j_2+i_1\leq 1$ that is not true.  It follows that $i_1+j_1=i_2+j_2$. Now, \eqref{100} gives $i_1=i_2$ and $j_1=j_2$. 
Therefore $\theta_{i_1} (\mathbb{N})\cap\theta_{i_2}(\mathbb{N})=\emptyset$, for all $i_1, i_2\in \mathbb{N}$ with $i_1\not= i_2$. Also, using definition \eqref{deftet} it is easy to see that maps $\theta_i$ are one-to-one for each $i\in\mathbb{N}$. Thus,
\begin{eqnarray}\label{theta_examp}
\Theta
:=\{\theta_k:\mathbb{N}\xrightarrow{1-1} \mathbb{N}\,|\;k\in I_0,\hspace{0.2cm}\theta_i (\mathbb{N})\cap\theta_j(\mathbb{N})=\emptyset,\hspace{0.2cm} i\neq j\}.
\end{eqnarray}

Let $T_1$ be an operator defined by
\begin{eqnarray}\label{form}
T_1=\summ_{i=1}^\infty\lambda_i P_{\theta_i},
\end{eqnarray}
where $\lambda=(\lambda_i)_{i\in \mathbb{N}}\in \ell^1(\mathbb{N})^+$ is an arbitrary fixed function.  The operator $T_1$   may be represented by an infinite matrix in the following way
$$
T_1=\left[\begin{matrix}
0 &  0  & 0&0 &0&\ldots  \\
\lambda_1 &  0  & 0&0 &0&\ldots  \\
0 &  \lambda_1  & 0&0 &0&\ldots  \\
\lambda_2 &  0  & 0&0 &0&\ldots  \\
   0&  0  & \lambda_1& 0&0 &\ldots  \\
0   &  \lambda_2    & 0 &  0  & 0&\ldots  \\
\lambda_3 &  0  & 0&0 &0&\ldots  \\
   0&  0  & 0 & \lambda_1& 0&\ldots  \\
      0&  0  & \lambda_2& 0&0 &\ldots  \\
0   &  \lambda_3    & 0 &  0  & 0&\ldots  \\
\lambda_4 &  0  & 0&0 &0&\ldots  \\
0   &  0 & 0&0 &\lambda_1&\ldots  \\
0   &  0  & 0&\lambda_2 &0&\ldots  \\
      0&  0  & \lambda_3& 0&0 &\ldots  \\
      0   &  \lambda_4    & 0 &  0  & 0&\ldots  \\
      \lambda_5 &  0  & 0&0 &0&\ldots  \\
\vdots & \vdots&\vdots&\vdots&\vdots& \ddots
\end{matrix}\right]
$$
Using \eqref{matrixrep} we obtain
$$Tf:=\left[0,\lambda_1f_1,\lambda_1f_2,\lambda_2f_1,\lambda_1f_3,\lambda_2f_2,\lambda_3f_1,\lambda_1f_4,\lambda_2f_3,\lambda_3f_2,\lambda_4f_1,\lambda_1f_5,\lambda_2f_4,\ldots\right]^T$$
for each $f\in\ell^p(\mathbb{N})$.
The operator $T_1$ is a bounded linear operator on $\ell^p(\mathbb{N})$, for all $p\in[1,\infty)$, by Theorem \ref{forallp}. 

Using $\ell^1(\mathbb{N})^+\subset \ell^p(\mathbb{N})^+$, we get that operator $T_1$ satisfies the statement $iii)$ of Corollary \ref{shape}. 
Also, using the matrix representation of $T_1$ we may conclude that $T_1$ satisfies statement $iv)$ of Corollary \ref{shape2}. Thus, $T_1$ preserves the submajorization relation on $\ell^p(\mathbb{N})^+$, for each $p\in[1,\infty)$.
Let
\begin{eqnarray}\label{form2}
T_2(f):=h\sum_{i\in \mathbb{N}}f(i),\;\;\forall f\in\ell^1(\mathbb{N})
\end{eqnarray}

where $h\in \ell^1(\mathbb{N})^+$ defined by 
\begin{equation*}\label{def_h}
h(j):=\left\{\begin{matrix}
a \geq 0, &  j=1,   \\
0, &  \text{otherwise}.   \\
\end{matrix}
\right.
\end{equation*}
Now, the operator $$ T:=T_1+T_2$$ 
 may be represented by an infinite matrix  
$$
T =\left[\begin{matrix}
a &  a  & a&a &a&\ldots  \\
\lambda_1 &  0  & 0&0 &0&\ldots  \\
0 &  \lambda_1  & 0&0 &0&\ldots  \\
\lambda_2 &  0  & 0&0 &0&\ldots  \\
   0&  0  & \lambda_1& 0&0 &\ldots  \\
0   &  \lambda_2    & 0 &  0  & 0&\ldots  \\
\lambda_3 &  0  & 0&0 &0&\ldots  \\
   0&  0  & 0 & \lambda_1& 0&\ldots  \\
      0&  0  & \lambda_2& 0&0 &\ldots  \\
0   &  \lambda_3    & 0 &  0  & 0&\ldots  \\
\lambda_4 &  0  & 0&0 &0&\ldots  \\
0   &  0 & 0&0 &\lambda_1&\ldots  \\
0   &  0  & 0&\lambda_2 &0&\ldots  \\
      0&  0  & \lambda_3& 0&0 &\ldots  \\
      0   &  \lambda_4    & 0 &  0  & 0&\ldots  \\
      \lambda_5 &  0  & 0&0 &0&\ldots  \\
\vdots & \vdots&\vdots&\vdots&\vdots& \ddots
\end{matrix}\right]
$$
and
$$Tf:=\left[a\summ_{i=1}^\infty f(i),\lambda_1f_1,\lambda_1f_2,\lambda_2f_1,\lambda_1f_3,\lambda_2f_2,\lambda_3f_1,\lambda_1f_4,\lambda_2f_3,\lambda_3f_2,\lambda_4f_1,\lambda_1f_5,\lambda_2f_4,\ldots\right]^T$$
for each $f\in\ell^1(\mathbb{N})$.
The operator $T $ is a bounded linear operator on $\ell^1(\mathbb{N})$ (by \cite[Theorem 3.1]{martin5}) because it satisfies \eqref{columnsup}. Now, $T$ preserves the submajorization relation on $\ell^1(\mathbb{N})^+$ by statement $iii)$ (or $iv)$) in Corollary \ref{shape2}.}
\end{example}

In the above example, we presented  linear preservers of submajorization on $\ell^1(\mathbb{N})^+$ which have only one "row" with mutually equal non-zero elements. In the next example we give  preservers where sets $I_1$ and $I_2$ in Corollary \ref{shape2} are both countable, that is, where there are countable "rows" which are a singleton.

\begin{example}\rm{
Let $\theta_i:\mathbb{N}\rightarrow\mathbb{N}$ be   maps defined by
\begin{eqnarray}\label{deftet2}
\theta_i(j)=i-1+\sum_{k=1}^{i+j-1}k,\;\;\;\;\;\forall  j\in \mathbb{N}
\end{eqnarray}
for any $i\in \mathbb{N}$. 
Similarly as in the previous example, we may provide that the family \eqref{theta_examp} contains one-to-one maps $\theta_i(j)$ with mutually disjoint images.

Fix a non-zero sequence $\mu=(\mu_i)_{i\in\mathbb{N}}\in\ell^1(\mathbb{N})^+$ and
suppose that $h\in\ell^1(\mathbb{N})^+$ is defined by
\begin{equation*}\label{def_h2}
h(j):=\left\{\begin{matrix}
\mu_i, &  j=\summ_{k=2}^{i+1} k,   \\
0, &  \text{otherwise}.   \\
\end{matrix}
\right.
\end{equation*}
In this way,  we define the operator $T:=T_1+T_2$, where  operators $T_1$ and $T_2$ are determined  as in \eqref{form} and  \eqref{form2}, respectively.

In order to provide that statement $iii)$ of Corollary \ref{shape2} is valid, suppose contrary that there exists $r\in \cup_{i\in \mathbb{N}}\theta_i(\mathbb{N})$ such that $\langle h,e_r\rangle>0$.
It follows that $r=  \theta_i(j)$ and $r=\summ_{k=2}^{n+1} k$ for some $i,j,n\in \mathbb{N}$. Hence, we get
$$ \summ_{k=2}^{n+1} k=r=i-1+\sum_{k=1}^{i+j-1}k$$ 
that is,
 $$ \summ_{k=1}^{n+1} k-\sum_{k=1}^{i+j-1}k=i\geq 1.$$ 
}
Therefore, $n+1> i+j-1$, and  so
$$i=\summ_{k=1}^{n+1} k-\sum_{k=1}^{i+j-1}k\geq i+j>i$$
which is a contradiction. Thus, $T\in  \mathcal{P}_{s}(\ell^1(I)^+)$, so the operator $T$ may be represented by an infinite matrix  in the following way
$$
T =\left[\begin{matrix}
\lambda_1 &  0  & 0&0 &0&\ldots  \\
\mu_1 &  \mu_1  & \mu_1&\mu_1 &\mu_1&\ldots  \\
0 &  \lambda_1  & 0&0 &0&\ldots  \\
\lambda_2 &  0  & 0&0 &0&\ldots  \\
\mu_2 &  \mu_2  & \mu_2&\mu_2 &\mu_2&\ldots  \\
   0&  0  & \lambda_1& 0&0 &\ldots  \\
0   &  \lambda_2    & 0 &  0  & 0&\ldots  \\
\lambda_3 &  0  & 0&0 &0&\ldots  \\
\mu_3 &  \mu_3  & \mu_3&\mu_3 &\mu_3&\ldots  \\
   0&  0  & 0 & \lambda_1& 0&\ldots  \\
      0&  0  & \lambda_2& 0&0 &\ldots  \\
0   &  \lambda_3    & 0 &  0  & 0&\ldots  \\
\lambda_4 &  0  & 0&0 &0&\ldots  \\
\mu_4 &  \mu_4  & \mu_4&\mu_4 &\mu_4&\ldots  \\
0   &  0 & 0&0 &\lambda_1&\ldots  \\
0   &  0  & 0&\lambda_2 &0&\ldots  \\
\vdots & \vdots&\vdots&\vdots&\vdots& \ddots
\end{matrix}\right].
$$
\end{example}


The submajorization relation $\prec_s$ may be defined on the whole Banach space $\ell^p(I)$ instead of the positive cone   $\ell^p(I)^+$. However, authors do not know whether preservers of such extension of submajorization will have the same shape as preservers of submajorization on $\ell^p(I)^+$. It seems that some of the main results such as Corollaries \ref{shape} and \ref{shape2}  will probably not be correct. As an   aditional argument for this, we notice that preservers of extended   weak majorization \cite{martin7}  on whole $\ell^1(I)$  do not have the same shape as preservers of weak majorization on $\ell^1(I)^+$ presented in \cite{martin4}.


\newpage

Address:

\bigskip
Martin Z. Ljubenovi\' c

Faculty of Mechanical Engineering, Department of Mathematics, University of Ni\v s, Aleksandra Medvedeva
14, 18000 Ni\v s, Serbia.
\bigskip

Dragan S. Raki\' c

Faculty of Mechanical Engineering, Department of Mathematics, University of Ni\v s, Aleksandra Medvedeva
14, 18000 Ni\v s, Serbia.
\bigskip

{\it E-mail}
\medskip

M. Z. Ljubenovi\'c: {\tt martinljubenovic@gmail.com}

D. S. Raki\' c:  {\tt rakic.dragan@gmail.com}

\end{document}